\documentclass[12pt]{amsart}
\usepackage{amsmath,amssymb,amscd,amsxtra}
\usepackage{latexsym}
\usepackage[dvips]{graphics,epsfig}

\newcommand{\R}{\mathbb{R}}

\newcommand{\Z}{\mathbb{Z}}
\newcommand{\M}[1]{\mathcal{M}^{#1}}
\renewcommand{\S}{\mathcal{S}}
\newcommand{\di}{\mathcal{S}^{'}}

\newcommand{\nm}[2]{\|{#1}\|_{#2}}

\newcommand{\bigparen}[1]{\bigl(#1\bigr)}

\newcommand{\biggparen}[1]{\biggl(#1\biggr)}
\newcommand{\ip}[2]{\langle#1,#2\rangle}

\newcommand{\zd}[1]{{\Z}^{{#1}d}}
\newcommand{\rd}{{\R}^{d}}

\newcommand{\modsp}{modulation space}
\newcommand{\wam}{W(\cF L^1, \ell ^\infty)}
\newcommand{\flone}{\cF L^1}
\newcommand{\cF}{\mathcal{F}}
\newcommand{\cS}{\mathcal{S}}
\newcommand{\intrd}{\int _{\rd }}

\newtheorem{thm}{Theorem}
\newtheorem{lemma}{Lemma}

\newtheorem{cor}{Corollary}
\theoremstyle{remark}
\newtheorem{rem}{Remark}
\theoremstyle{definition}
\newtheorem{deft}{Definition}

\begin{document}

\title{Unimodular  Fourier multipliers for modulation spaces}

\author{\'Arp\'ad B\'enyi}

\address{\'Arp\'ad B\'enyi\\
Department of Mathematics\\
516 High Street\\
Western Washington University\\
Bellingham, WA 98225, USA }

\email{arpad.benyi@wwu.edu}

\author{Karlheinz Gr\"ochenig}

\address{Karlheinz Gr\"ochenig\\
Faculty of Mathematics\\
Universit\"at Wien\\
1090 Vienna, Austria }

\email{karlheinz.groechenig@univie.ac.at}

\author{Kasso A.~Okoudjou$^*$}

\address{Kasso A.~Okoudjou\\
Department of Mathematics\\
University of Maryland\\
 College Park, MD 20742, USA}

\email{kasso@math.umd.edu}

\author{Luke G. Rogers}

\address{Luke G.~Rogers\\
Department of Mathematics\\
Malott Hall\\
Cornell University\\
Ithaca, NY 14853, USA}

\email{luke@math.cornell.edu}

\thanks{$^*$ Research partially supported by an Erwin Schr{\"o}dinger
  Junior Fellowship and by NSF grant DMS-0139261. K.~G.~was supported
  by the Marie-Curie Excellence Grant MEXT-CT 2004-517154.}
 \subjclass[2000]{Primary 42B15;
Secondary 42B35, 47G30}

\date{\today}

\keywords{Fourier multiplier, modulation space, short-time Fourier transform, Schr\"odinger
equation, wave equation, conservation of energy}

\begin{abstract}
We investigate the boundedness of unimodular  Fourier multipliers on modulation spaces.
Surprisingly, the  multipliers with general symbol $e^{i|\xi|^\alpha }$, where $\alpha\in[0, 2]$,
are bounded on all modulation spaces, but, in general, fail to be bounded on the usual
$L^p$-spaces. As a consequence, the phase-space concentration of the solutions to the free
Schr\"odinger and wave equations are preserved. As a byproduct, we also obtain boundedness results
on modulation spaces for singular multipliers $|\xi|^{-\delta} \sin (|\xi|^\alpha)$ for  $ 0 \leq
\delta \leq \alpha $.
\end{abstract}

\maketitle \maketitle \pagestyle{myheadings} \thispagestyle{plain} \markboth{\'A. B\'enyi, K.
Gr\"ochenig, K. A. Okoudjou and L. G. Rogers} {Unimodular Fourier Multipliers}

\section{Introduction and Motivation} A Fourier multiplier is a linear
operator $H_\sigma $  whose action on  a test  function $f$ on $\rd$  is formally defined by
\begin{equation}\label{fourmulti}
H_{\sigma} f(x) = \int_{\rd{}}\sigma(\xi)\, \hat{f}(\xi)\, e^{2 \pi i \xi \cdot x}\, d\xi \, .
 \end{equation}
The function $\sigma $ is called the  symbol or multiplier. Using the inverse Fourier transform,
one can also rewrite the operator as a convolution operator
$$H_{\sigma}
(x) = \check{\sigma} \ast f (x)$$ where $\check{\sigma}$ is the (distributional)  inverse Fourier
transform of $\sigma$.

Fourier multipliers arise naturally in the  formal solution of linear PDEs with constant
coefficients and in the convergence of Fourier series. For this reason, a fundamental problem in
the study of Fourier multipliers is to relate the boundedness properties of $H_\sigma$ on certain
function spaces to properties of the symbol $\sigma$. On $L^{2}$, $L^{1}$ and $L^{\infty}$ this is
relatively straightforward, while the full resolution of this problem for general $L^p$-spaces is
an analytic gem known as the H\"ormander-Mihlin multiplier theorem~\cite{Ho63,mi56}. A detailed
exposition of the theory of Fourier multipliers may be found in \cite{Gra03, Ste70}.

In this paper we study unimodular Fourier multipliers, in particular the multipliers with symbol
$e^{it|\xi |^\alpha }$ for $t \in \mathbb{R}$ and $\alpha \in [0,2]$. To understand the problem,
consider first a general unimodular multiplier $\sigma (\xi )= e^{im(\xi )}$ for some real-valued
function $m$. If $|\partial ^\alpha m (\xi )| \leq C_\alpha |\xi |^{-\alpha} $ for all $\xi \neq 0$
and $|\alpha | \leq [d/2]+1$, then $e^{im}$ is bounded on $L^p$ by Mihlin's condition,
see~\cite[p.~367]{Gra03}. The case of the  unimodular function $e^{i|\xi |^\alpha }$ is more
complicated.  In addition to the singularity of the derivatives at the origin, the  multiplier has
large  oscillations at infinity and thus possesses large derivatives. If $\alpha > 1$, this
excludes the application of the  multiplier theorems of H\"ormander-Mihlin and their variations.
Indeed, the multiplier may be unbounded. Specifically, if  $\alpha >1$,  then the operator
$H_{e^{i|\xi |^\alpha}}$ is  bounded on $L^{p}(\rd{})$ if and only if  $p=2$ as a consequence
of~\cite{Ho63, LeOl}.

The cases $\alpha =1$ and $\alpha =2$ are particularly interesting and have been studied
intensively in PDE, because they occur in the time evolution of the wave equation ($\alpha = 1$)
and the free Schr\"odinger operator ($\alpha=2$). The unboundedness of the multiplier on general
$L^p$ means that $L^p$-properties of the initial conditions are not preserved by the time
evolution.   There is an extensive literature about  estimates on the wave operator; see, for
example, \cite{Be82, Lit63, Per, Str77}.   Little is known for $\alpha \in (0, 1)$, and many
results are concerned almost exclusively with appropriate corrections of the symbol $e^{i|\xi
|^\alpha }$ with functions that are essentially smooth away from the origin, \cite{Hir, Per, Wa65}.

In view of the unboundedness of the multiplier $e^{i|\xi |^\alpha }$ on $L^p$ it is natural to
question whether the $L^p$-spaces are really the appropriate function spaces for the study and
understanding of these operators. We suggest that they are not, and propose that the so-called
\emph{modulation spaces} are a good alternative class for the study of unimodular Fourier
multipliers.

To define the modulation spaces we fix a non-zero Schwartz function $g$ and consider the short-time
Fourier transform $V_g f$ of a function $f$ with respect to $g$
$$
    V_g f(x, \omega)
    =\ip{f}{M_\omega T_x g} =\int_{\rd{}} e^{-2\pi i \omega \cdot  t} \,
        \overline{g(t-x)} \, f(t) \, dt \, .
$$
The modulation space $\M{p,q}$ is the closure of the Schwartz class with respect to the norm
$$
\nm{f}{\M{p,q}}=\biggparen{\int_{\rd{}}\biggparen{\int_{\rd{}}|V_{g}f(x,\omega)|^p\,
dx}^{q/p}\, d\omega}^{1/q} 
$$
(with appropriate modifications  when $p=\infty $ or $q=\infty $). A priori it may seem that this
norm depends on $g$, so it is worth noting that different choices of $g$ give equivalent norms.

Since their introduction by  Feichtinger~\cite{Fei83}, modulation spaces have become canonical for
both time-frequency and phase-space analysis.  Their many applications are surveyed in
\cite{fei06}; for the special case of $\M{\infty ,1}$, which is sometimes called the Sj\"ostrand
class, see also~\cite{grhe, gro06, Sj94, toft} and the references therein.

The reason for the ubiquity of modulation spaces is essentially that $V_gf$ is a local version of
the Fourier transform. In the terminology of physics, if $x$ is the position and $\omega$ the
momentum of a physical state, then $V_gf (x,\omega)$ is a measure of the amplitude of a state $f$
at the point $(x,\omega )$ in phase space. The modulation space norm $\|\cdot \|_{\M{p,q}}$ can
then be understood as a measure for the phase space concentration of $f$.  In this interpretation,
boundedness of a Fourier multiplier on modulation spaces expresses the conservation of phase-space
properties, which is the natural extension of the energy conservation corresponding to the obvious
$L^2$-boundedness.

An abstract characterization  of all Fourier multipliers on modulation spaces was obtained in
\cite{FeiNa} (see also Theorem \ref{genemulthm} below), however this characterization requires a
deep understanding of Fourier multipliers on $L^p$-spaces and it is often very difficult to check
that the given conditions are satisfied for a given symbol.  One exception is the case where the
multiplier has sufficiently many bounded derivatives (\cite{FeiNa} Theorem 20) where the
H\"{o}rmander-Mihlin theorem can be applied locally. We note that this result is not applicable to
the multipliers $e^{it|\xi|^{\alpha}}$.

Our main result is that the unimodular multipliers discussed above are bounded on all modulation
spaces.

\begin{thm}\label{maintheorem}
If  $\alpha \in [0,2]$, then the Fourier multiplier $H_{\sigma_{\alpha}}$  with
$\sigma_{\alpha}(\xi) = e^{i |\xi|^{\alpha}}$is   bounded from $\M{p,q}(\rd{})$ into
$\M{p,q}(\rd{})$  for all $1\leq p, q \leq \infty$ and in any dimension  $d\geq 1$.
\end{thm}

 As a consequence, we will
show that the Cauchy problems for the free  Schr\"odinger equation  and the wave equation with
initial data in a modulation space  satisfy  an analog of the principle of conservation of energy.
It is worth noting that modulation spaces were recently rediscovered and many of their properties
reproved in~\cite{wang06} in a study of the non-linear Schr\"odinger equation and the
Ginzburg-Landau equation. In particular, Corollary~\ref{humm}(a) was obtained in~\cite{wang06} by a
different method than that used here.

\begin{cor}\label{humm}

(a) Let $u(x,t)$ be the solution of the free Schr\"odinger equation
  $iu_t = \Delta _x u$ and $u(x,0) = f(x)$. If $f\in \M{p,q}$, then
  $u(\cdot, t)\in \M{p,q}$ for all $t>0$.

(b)  Let $u(x,t)$ be the solution of the wave  equation
  $iu_{tt} = \Delta _x u$ and $u(x,0) = f(x)$, $u_t(x,0) = g(x)$ . If
  $f, g\in \M{p,q}$, then
  $u(\cdot,t)\in \M{p,q}$ for all $t>0$.
\end{cor}

 To put it more succinctly, we may say that \emph{the phase-space
   concentration of an initial state  is preserved under the time
   evolution of the free Schr\"odinger equation and the wave equation.}
This result  is in striking contrast to the behavior of these Cauchy problems with initial data in
$L^p$ spaces \cite{Ho63, Lit63, Str70a, Str70b, Str77}.

 Our paper is organized as follows. In
Section $2$ we set up the notation and define the modulation spaces and some amalgam spaces needed
for the multiplier theory.  Section $3$ is devoted to the abstract characterization of Fourier
multipliers of modulation spaces given in \cite{FeiNa}. Furthermore, several sufficient conditions
for a Fourier multiplier to be bounded on the modulation spaces are provided. These conditions are
then used to prove our main results, which are stated and proved in Section 4. Finally, Section $5$
deals with the applications of our results to the analysis of the Cauchy problems associated to the
Schr\"odinger and wave equations.

\section{The Short-Time Fourier Transform and Associated Function
  Spaces}

\subsection{General notation}

For $x, \omega \in \rd{}$ the translation and modulation operators acting on a function $f$ defined
over $\rd{}$ are given respectively by  $$T_x f(t) = f(t-x) \qquad{\rm and}\qquad M_\omega f(t) =
e^{2\pi i \omega \cdot  t} \, f(t), $$ where $t\in \rd{}$. The Schwartz class of test functions
will be denoted by  $\S =\S(\rd{})$, its dual is the space  of  tempered distributions  $\di=
\di(\rd{})$  on $\rd{}$. The Fourier transform of $f \in \S$ is given by  $$\hat{f}(\omega)
=\int_{\rd{}} f(t) \, e^{-2\pi i t \cdot \omega} \, dt, \, \omega \in \rd{},$$  which is an
isomorphism of the Schwartz space $\S $ onto itself that extends to the tempered distributions by
duality. The inverse Fourier transform is given explicitly  by $\check{f} (x) = \intrd f(\omega )
e^{2\pi i x\cdot\omega } \, d\omega $, and we have $(\hat{f})\check{} = f$. The inner product of
two functions $f,g\in L^2$ is $\ip{f}{g}=\int_{\rd{}}f(t)\overline{g(t)} \,dt$, and its extension
to $\di\times \S$ will be also denoted by $\ip{\cdot}{\cdot}$.

A key object in time-frequency analysis is the short-time Fourier transform (STFT), which in a
sense is a ``local'' Fourier transform that has the advantage of displaying the frequency content
of any function during various time intervals. More precisely, the STFT of a tempered distribution
$f \in \di$ with respect to a window $g \neq 0 \in \S$ is
\begin{equation}
  \label{eq:ch11}
  V_g f(x, \omega)
=\ip{f}{M_\omega T_x g} =\int_{\rd{}} e^{-2\pi i \omega \cdot  t} \, \overline{g(t-x)} \, f(t) \,
dt \, ,
\end{equation}
(where the integral version exists only for functions of polynomial growth). We will consistently
use the following equivalent forms for the STFT
\begin{equation}\label{eq:STFTequivalences}
    V_{g}f(x, \omega) = \widehat{f\cdot T_{x}\overline{g}}(\omega) = e^{-2\pi i x \cdot
    \omega}V_{\hat{g}}\hat{f}(\omega, -x).
    \end{equation}
If $g\in \cS $ and $f\in \cS '$, then $V_gf $ is a continuous function of polynomial
growth~\cite{Gr01}. In a less obvious way, the STFT can   be defined even   when both $f\in \S '
(\rd{} )$ and $g\in \S ' (\rd{} )$ \cite[Prop.~1.42]{Fol89}.

The time-frequency content of a tempered distribution can be quantified by imposing a mixed-norm on
its STFT. Throughout the paper, we let $L^{p,q}=L^{p,q}(\rd{}\times \rd{})$ be the spaces of
measurable functions $f(x,\omega)$ for which the mixed norm
$$\nm{f}{L^{p,q}}=\biggparen{\int_{\rd{}}\biggparen{\int_{\rd{}}|f(x, \omega)|^p\,
dx}^{q/p}d\omega}^{1/q}$$ is  finite.  If $p=q$, we have $L^{p,p}=L^p$, the usual Lebesgue spaces.

We use the notation $u\lesssim v$ to denote $u\leq cv$, for a universal (independent of $u$ and
$v$) positive constant $c$. Similarly, we use the notation $ u \asymp v$ to denote $c u \leq v \leq
C u$, for some universal positive constants $c, C$.

\subsection{Modulation spaces}
\begin{deft}
Given $1\leq p, q \leq \infty$, and given a non-zero  window function $g \in \S$, the modulation
space $\M{p,q}=\M{p,q}(\rd{})$ is the space of all distributions $f \in \di$ for which the
following norm is finite:
\begin{equation}
\label{mod} \nm{f}{\M{p,q}}=\biggparen{\int_{\rd{}}\biggparen{\int_{\rd{}}|V_{g}f(x,\omega)|^p\,
dx}^{q/p}\, d\omega}^{1/q} = \nm{V_{g}f}{L^{p,q}},
\end{equation}
with the usual modifications if $p$ and/or $q$ are infinite. When $p=q$, we will write $\M{p}$ for
the
 modulation
space $\M{p,p}$.
\end{deft}

 This definition is independent of the choice of the window $g$ in
the sense of equivalent norms. Moreover, if $1\leq p, q <\infty$, then $\M{1}$ is densely embedded
into $\M{p,q}$, as is the Schwartz class $\S$. If $1\leq p, q < \infty$, then the dual of $\M{p,q}$
is $\M{p', q'}$, where $1\leq p, q < \infty$ and $\tfrac{1}{p}+\tfrac{1}{p'}=
\tfrac{1}{q}+\tfrac{1}{q'}=1$. We refer to \cite{Fei83, Gr01} and the references therein for the
precise  details and the rich theory of modulation spaces.

Note that an application of Plancherel's  theorem yields $\M{2}=L^2$. However, it can be shown that
for $p, q \neq 2$, $\M{p,q}$ does not coincide with any Lebesgue space. Instead, one may use the
embeddings $\M{p}\subset L^p  \subset \M{p,p'}$, if $1\leq p\leq 2 $ and $\M{p,p'}\subset L^p
\subseteq \M{p}$ for $p\geq 2$. We will also use the fact that modulation spaces are invariant
under dilation, i.e., if $f \in \M{p,q}$, the $f(t \cdot ) \in \M{p,q}$ for every $t >0$, see for
instance ~\cite[Ch.~9]{Gr01}.

\subsection{Wiener Amalgam Spaces}

If we reverse the order of integration in~\eqref{mod}, then we obtain the Wiener amalgam
spaces~\cite{fei90}. Let us write $\flone$ for the space of all Fourier transforms of $L^1$, that
is
\begin{equation*}
  \label{eq:c1}
  \flone = \{ f \in L^\infty (\rd ): \hat{f} \in L^1(\rd ) \},
\end{equation*}
with norm $\|f \|_{\flone } = \|\hat{f} \|_{L^1}$.  This allows us to define an amalgam space that
will be used frequently.
\begin{deft}\label{amalgam}
Fix  $g \in \cS (\rd ), g \neq 0$. Then  the space $\wam $ consists of all functions $\sigma \in
L^\infty (\rd )$ for which
\begin{eqnarray}
  \| \sigma \|_{\wam } &=& \sup _{x\in \rd } \intrd  |V_g \sigma
  (x,\omega )|\, d\omega \notag  \\
&=& \sup _x \| \sigma \, T_x g \|_{\flone } <\infty .       \label{eq:c2}
\end{eqnarray}
\end{deft}

Roughly speaking, $\sigma $ is in $\wam $ if $\sigma $ is locally in the Fourier algebra $\flone $
with the local norms being uniformly bounded. As in the case of   modulation spaces, the definition
of $\wam$  is independent of the test function $g \in \cS $. We refer to \cite{fei90} and the
references therein for more details on these Wiener amalgam-type spaces.

For the abstract multiplier theorem of Feichtinger-Narimani (Theorem \ref{genemulthm} below) we
need one other Wiener amalgam space.
\begin{deft}
For fixed  $g \in \cS (\rd ), g \neq 0$ the space $W(M_{\mathcal{F}}(L^{p}), \ell^\infty)$ consists
of all tempered distributions $\sigma$ for which
\begin{equation}\label{wfmlp}
\nm{\sigma}{W(M_{\mathcal{F}}(L^{p}), \ell^\infty)} = \sup_{x \in
  \rd{}} \nm{H_{\sigma \, T_{x} g}}{L^p \to L^p} <\infty .
\end{equation}
 In particular, each $\sigma \in W(M_{\mathcal{F}}(L^{p}), \ell^\infty)$  coincides
locally with a
 Fourier multiplier on $L^p$.
\end{deft}

\section{Abstract Multiplier Theorems on Modulation Spaces}

The results of this paper were in part inspired by earlier work of three of the authors and Loukas
Grafakos \cite{BGGO}, in which an analogue of the classical Marcinkiewicz multiplier theorem was
proven in the modulation space context. By an easy tensor product argument, the proof of Theorem 1
of \cite{BGGO} (see also Corollary 19 of \cite{FeiNa}) may be extended to show the following.

\begin{thm}\label{hidbggo}For any $b=(b_1, b_2, \hdots, b_d) \in (\R_{+})^{d}$ with all $b_j>0$,
let $Q_{b}=\prod_{j=1}^{d} (0, b_j)$.  For a bounded sequence
 ${\bf c}=(c_n)_{n \in \zd{}} \in \ell^{\infty}(\zd{}),$ define the function
 \begin{equation*}
 \sigma_{b, {\bf c}} = \sum_{n \in \zd{}}c_{n} \chi_{n+ Q_{b}},
\end{equation*}
where $\chi_{E}$ is the indicator function of the set $E$.  Then the operators $H_{\sigma_{b, {\bf
c}}}$ are bounded from $\M{p,q}$ into $\M{p,q}$, $1 <p < \infty$, $1\leq q \leq \infty$, with a
norm estimate
$$\nm{H_{\sigma_{b, {\bf c}}}f}{\M{p,q}} \leq C(b,p,q) \nm{\bf
  c}{\ell^{\infty}} \nm{f}{\M{p,q}}\, .$$
\end{thm}

This theorem gives a concrete example of a multiplier that is not bounded on $L^p$ for $p\neq 2$
except in trivial cases.  To see that it is a version of the usual Marcinkiewicz theorem, the
reader should note that the proof of the one-dimensional case in \cite{BGGO} uses only that
$\sigma$ is bounded and of bounded variation, and that there is a uniform bound
$\int_{nb}^{(n+1)b}|d\sigma|\leq C$, $n\in\mathbb{Z}$, on the variation on $b$-length intervals.
Careful examination of the proof also shows that the boundedness on modulation spaces is related to
the localization of the multiplier in time and frequency.  It was this idea that led to Theorem
\ref{lemma:c2d} below.

A key feature of modulation spaces that permits boundedness of a class of multipliers larger than
that for $L^{p}$ is the fact that multipliers for different locations act approximately
independently. Feichtinger and Narimani \cite{FeiNa} made this idea precise to give an abstract
characterization of all Fourier multipliers on modulation spaces.

\begin{thm}[\protect{\cite{FeiNa} Theorem 17(1)}]\label{genemulthm}
A multiplier is bounded on $\M{p,q}$ if and only if  $\sigma \in W(M_{\mathcal{F}}(L^{p}),
\ell^\infty)$.
\end{thm}

Since this characterization rests on understanding all Fourier multipliers on $L^p$ it is difficult
to apply to a concrete function. Note, however, that Theorem 20 of \cite{FeiNa} shows that a
multiplier with $\lfloor d/2 \rfloor +1$  bounded derivatives on $\mathbb{R}^{d}$ is bounded on the
modulation spaces.  This result could be used to prove that part of Theorem \ref{maintheorem} which
is established in Theorem \ref{lemma:c2d} in the case $\alpha<1$, but is not useful in general for
the multipliers $e^{it|\xi|^{\alpha}}$ for $\alpha\in [0,2]$ because there is a singularity at the
origin and there are high frequencies (unbounded derivatives) far from the origin when $\alpha>1$.

As an alternative to the abstractness of Theorem \ref{genemulthm}, we offer several sufficient
conditions that are easier to verify in practice.  In particular we will see in Theorem
\ref{lemma:c2d} that they are valid when the multiplier is well localized in time-frequency space.
We remark that the conditions in the next result are far from being necessary.

\begin{lemma}\label{suffmulthm}
The Fourier multiplier  $H_\sigma $ is bounded on all \modsp s $\M{p,q}(\rd{})$ for $d\geq 1$ and
$1\leq p,q \leq \infty $ under each of the following conditions:

(i) $\sigma \in \wam $.

(ii)  $\sigma \in \M{\infty ,1}$.

(iii) $\sigma \in \flone $.
\end{lemma}

\begin{proof}
(i) Since $\flone \subseteq M_{\mathcal{F}}(L^p)$, we have $\wam \subset W( M_{\mathcal{F}}(L^p),
\ell ^\infty )$ and the claim follows from Theorem~\ref{genemulthm}. However, to keep the
presentation self-contained, we sketch a simple and direct proof that is based on the convolution
relations in~\cite{CG03}.

Let $g^*(x) = \overline{g(-x)}$. Then we can write the modulus of the STFT as $|V_gf(x,\omega )| =
|f \ast M_\omega g^*)(x)|$ and the modulation space norm as
$$
\|f\|_{\M{p,q}} = \Big( \int _{\rd } \|f \ast M_\omega g^* \|_{L^p}^q \, d\omega \Big)^{1/q} \, .
$$
Likewise
\begin{eqnarray*}
\|\sigma \|_{\wam } &\asymp &\sup _\omega \int _{\rd } |\langle \sigma,
M_{-x} T_\omega \hat{g} \rangle | \, dx \\
&=&    \sup _\omega \int _{\rd } |\langle \check{\sigma}, T_{x} M_\omega g \rangle | \, dx = \sup
_\omega \|\check{\sigma } \ast M_\omega g ^* \|_{L^1} \, .
\end{eqnarray*}
Now we choose a window that factors as $g_1 = g \ast g$ for some $g\in \cS $ and observe that
$M_\omega g_1 = M_\omega g \ast M_\omega g$. The boundedness of $H_\sigma $ for $\sigma \in \wam $
now follows from the following chain of inequalities, where we use repeatedly that
 the modulation space norm is independent of the window.
\begin{eqnarray*}
  \|H_\sigma f \|_{\M{p,q}} ^q &=& \int _{\rd } \| (H_\sigma f) \ast
  M_\omega g_1^* \|_{L^p} ^{q} \, d\omega \\
&=&  \int _{\rd } \| \check{\sigma } \ast f \ast
  M_\omega g^*\ast M_\omega g^* \|_{L^p} ^{q} \, d\omega \\
&\leq & \int _{\rd } \| \check{\sigma } \ast
  M_\omega g^* \|_{L^1} ^{q}   \|f \ast M_\omega g^* \|_{L^p}^q \,
  d\omega \\
&\leq &\sup _\omega \| \check{\sigma } \ast
  M_\omega g^* \|_{L^1} ^{q} \, \int _{\rd }    \|f \ast M_\omega g^*
  \|_{L^p}^q \,   d\omega \\
&\asymp & \|\sigma \|_{\wam }^{q} \, \|f \|_{\M{p,q}}^{q} \, .
\end{eqnarray*}

Statements (ii) and (iii)  now  follow from the embeddings $\flone \subset
  \M{\infty, 1} \subset \wam $.
For the first embedding, if we let $ g\in \S, g\neq 0$,  and $f \in \mathcal{F}L^1$, then by
\eqref{eq:STFTequivalences}
 \begin{align*}
    \int_{\rd{}}\sup_{x\in \rd{}}|V_{g}f(x, \omega)|\, d\omega
    &=\int_{\rd{}}\sup_{x\in\rd{}}|V_{\hat{g}}\hat{f}(\omega, -x)|\, d\omega\\
    &= \int_{\rd{}}\sup_{x\in \rd{}}|\widehat{\hat{f}\, \overline{T_{\omega}\hat{g}}}(-x)|\, d\omega\\
    &\leq \int_{\rd{}}\int_{\rd{}}|\hat{f}(y)\overline{T_{\omega}\hat{g}(y)}|\,dy\, d\omega\\
    &\leq \int_{\rd{}}|\hat{f}|\ast|\hat{g}|(\omega)\, d\omega\\
    &\leq \nm{\hat{f}}{L^{1}}\, \nm{\hat{g}}{L^{1}}.
\end{align*}
The second embedding is even easier:
\begin{eqnarray*}
\|f \|_{\wam }&=&\sup _{x\in \rd } \intrd  |V_g f (x,\omega
  )| \, d\omega\\
&\leq & \intrd \sup _{x\in \rd } |V_g f (x,\omega
  )| \, d\omega =   \|f\|_{\M{\infty,1}} \, .
\end{eqnarray*}
\end{proof}

\section{Unimodular Functions as Fourier Multipliers}

The primary goal of this section is to prove Theorem \ref{maintheorem}, showing that the
multipliers $e^{i|\xi|^{\alpha}}$ are bounded on all modulation spaces for $\alpha\in[0,2]$. Except
in a few special cases these multipliers  fail to be bounded on $L^p, p\neq 2$.

As earlier noted, there are two main obstacles to overcome, namely the singularity at $\xi = 0$
(except when $\alpha $ is an even integer) and large oscillations at infinity (for $\alpha>1$)
which give large derivatives  and preclude application of the  multiplier theorems of
H\"ormander-Mihlin and their variations. For clarity we will treat the singularity at the origin
separately from the oscillations at infinity.

Our methods  will also establish a stronger time-frequency property  of the multipliers when
$\alpha \in [0, 1]$, and will show that the singular multipliers $\Im
(e^{i|\xi|^\alpha}|\xi|^{-\delta})= |\xi |^{-\delta}\sin |\xi |^\alpha $, $\alpha \geq \delta\geq
0$, are bounded on the modulation spaces.

In what follows, $\chi \in \mathcal{C}^{\infty}_{c}(\rd{})$ will denote a test function such that
\begin{equation}\label{bump}
    \chi(\xi) = \begin{cases}
        0 \quad&\text{if $|\xi| \geq 2$}\\
        1 &\text{if $|\xi| \leq 1$}\\
        0 \leq \chi(\xi) \leq 1 &\text{if $1\leq | \xi | \leq 2$.}
        \end{cases}
\end{equation}

\subsection{The Singularity at the Origin}

For the purposes of establishing boundedness, the relevant feature of the singularity at the origin
is its homogeneity.
\begin{thm}\label{lem:ch4}
Assume that $\mu \in \mathcal{C}^{d+1}(\rd{}\setminus \{0\})$ is  homogeneous of order $\alpha >0$,
that is, $\mu (s \xi ) = s^\alpha \mu (\xi )$ for all $s>0$ and all $\xi \neq 0$. Then $e^{i \mu}
\chi \in \flone $, and consequently $H_{e^{i\mu }\chi }$ is bounded on all \modsp s
$\M{p,q}(\R^d)$, as well as on all Lebesgue spaces $L^p(\R^d)$ for $1\leq p,q \leq \infty $.
\end{thm}

\begin{proof}
We expand $e^{i\mu(\xi)}\chi$ as
\begin{equation*}
    e^{i\mu(\xi)} \chi(\xi)
    = \sum_{k=0}^{\infty} \frac{i^k}{k!}\mu(\xi)^{k}\chi(\xi)
\end{equation*}
and will show that $\phi_{k} = \mu^{k}\chi \in \mathcal{F}L^1$ for all  $k\in \mathbb{N}$ with norm
estimates sufficient to ensure convergence.  To do so, define $\psi(\xi) = \chi(\xi/2) -
\chi(\xi)$, so that $ supp \, (\psi) \subset \{\xi \in \rd{}: 1 \leq |\xi| \leq 4\}$ and
$\sum_{j=1}^{\infty} \psi(2^{j}\xi) = \chi(\xi)$ for all $\xi \neq 0$. Using this and the
homogeneity of
 $\mu$ we further decompose $\phi_k$ as
\begin{equation}\label{eq:estforphik}
    \phi_{k}(\xi)
    = \sum_{j=1}^{\infty}\mu(\xi)^{k}\psi(2^{j}\xi)
    =\sum_{j=1}^{\infty}2^{-kj\alpha}\,\mu(2^{j}\xi)^{k}\psi(2^{j}\xi)
    =\sum_{j=1}^{\infty}2^{-kj\alpha}\,\psi_{k}(2^{j}\xi)
\end{equation}
where $\psi_{k} = \mu^{k}\psi \in\mathcal{C}^{d+1}$ and has compact support. The Fourier transform
of $\psi_{k}(2^{j}\xi)$ is $2^{-jd}\hat{\psi}_{k}(2^{-j}\xi)$, so that $\|\psi _k (2^j
\cdot)\|_{\cF L^1}=\|\hat{\psi _k}\|_{L^1}=\|\psi _k \|_{\cF L^1}$ is independent of $j$. We
estimate it in two parts
$$
\nm{\hat{\psi}_{k}}{L^{1}}  =\int_{|\omega| \leq 1} |\hat{\psi_{k}}(\omega)|\, d\omega +
\int_{|\omega|\geq 1} |\hat{\psi_{k}}(\omega)|\, d\omega = A +B.$$

To estimate $A$, observe that homogeneity of $\mu$ ensures there is a constant $C_{0}$ such that
$\mu^k$ is bounded by $C_{0}^{k}4^{k\alpha}$ on the support of $\psi$, and compute
\begin{equation}\label{eq:estforA}
    A
    \leq
    v_{d} \nm{\hat{\psi_{k}}}{L^{\infty}}
    \leq v_d \sup _{|\xi | \leq 4} |\mu (\xi )|^k \, \|\psi \|_{L^1}
    \leq
    v_{d} C_0^k 4^{k\alpha}\nm{\psi}{L^{1}}
    \lesssim C_0^k 4^{k\alpha}
\end{equation}
where $v_d$ is the volume of the unit ball in $\rd{}$.

The estimate for $B$ is a little more involved.  Since $\widehat{\psi_{k}}(\omega) = (2\pi i
\omega)^{-\beta} \widehat{\partial^{\beta}\psi_{k}}(\omega)$ for all $|\beta| \leq d+1$ we may use
the pointwise estimate
\begin{eqnarray*}
|\widehat{\psi_{k}}(\omega)| &\leq & \min _{|\beta | \leq d+1} |(2\pi i
\omega)^{-\beta}| \, |\widehat{\partial^{\beta}\psi_{k}}(\omega)| \\
&\lesssim & \max _{|\beta | \leq d+1} \nm{\widehat{\partial^{\beta}\psi_{k}}}{L^\infty } \min
_{|\beta | \leq d+1} | \omega^{-\beta}| \, .
  \end{eqnarray*}
Then we have
 \begin{equation}
    B
    \leq \max _{|\beta |\leq d+1} \nm{\widehat{\partial^{\beta}\psi_{k}}}{L^{\infty}}\int_{|\omega|\geq 1}\frac{1}{|\omega^{\beta}|}\, d\omega
    \leq \max _{|\beta |\leq d+1}
    \nm{\partial^{\beta}\psi_{k}}{L^{1}}\int_{|\omega|\geq 1} \min
    _{|\beta |\leq d+1} \frac{1}{|\omega^{\beta}|}\, d\omega.
  \label{estb1}
  \end{equation}
By Leibniz's rule, the derivative $\partial^\beta (\mu ^k \psi )$ is a sum of
$\binom{k+|\beta|}{k}$ terms of the form $(\partial ^{\gamma _{k+1}}\psi)\prod _{j=1}^k \partial
^{\gamma _j} \mu$ with $\sum _{j=1}^{k+1} \gamma _j = \beta $.  Each of those involving $\mu$ may
be estimated using homogeneity, since
\begin{equation*}
    \partial^{\gamma_{j}}\mu(\xi)=s^{-\alpha}s^{|\gamma_{j}|}\big(\partial^{\gamma_{j}}\mu\big)(s\xi)
    \qquad \forall s>0, \xi \neq 0 \, .
    \end{equation*}
In particular we see that when $|\gamma_{j}|\leq d+1$ and $1\leq |\xi|\leq 4$
\begin{equation*}
    |\partial^{\gamma_{j}}\mu^{k}(\xi)|
    \leq 4^{\alpha-|\gamma_{j}|} \sup_{|\xi|=1} |\partial^{\gamma_{j}}\mu(\xi)|
    \leq C_1 4^{\alpha-|\gamma_{j}|}
    \end{equation*}
where $C_1$ is a bound for the first $d+1$ partial derivatives of $\mu$ on the unit sphere.  Such a
bound exists because $\mu\in \mathcal{C}^{d+1}(\mathbb{R}^{d}\setminus{0})$.  The derivatives of
$\psi$ are clearly bounded, so each term of $\partial^\beta (\mu ^k \psi )$ satisfies $|\partial
^{\gamma _{k+1}}\psi|\prod _{j=1}^k |\partial ^{\gamma _j} \mu| \lesssim 4^{k\alpha}$.  Crudely
estimating the number of these by $\binom{k+|\beta|}{k}\leq C_{2}^k$ for a sufficiently large
$C_2=C_{2}(d)$ we arrive at the bound $\max_{|\beta|\leq d+1}
\|\partial^{\beta}\psi_{k}\|_{L^{1}}\lesssim 4^{k\alpha}C_{2}^k$.

It is not difficult to show that $\min_{|\beta|\leq d+1} |\omega^{\beta}|^{-1}$ is integrable
outside the unit sphere (see, e.g., \cite[pp. 321]{Gr01}), so it follows from \eqref{estb1} that
$B\lesssim C^k 4^{k\alpha}$.  Combining this with \eqref{eq:estforA} we obtain
\begin{equation*}
    \nm{\widehat{\psi}_{k}}{L^{1}}\lesssim C^k 4^{k\alpha}
    \end{equation*}
and therefore from \eqref{eq:estforphik},
$$\nm{\phi_k}{\flone}=\nm{\hat{\phi}_{k}}{L^{1}}\leq \sum _{j=1} ^\infty 2^{-kj\alpha }
\|\widehat{\psi _k} \|_{L^1}  \lesssim C^{k} \frac{2^{k\alpha}}{1-2^{-k\alpha}}.$$ Consequently,
$$\|e^{i\mu } \chi \|_{\flone } \leq \sum _{k=0}^\infty
\frac{1}{k!} \, \|\phi _k \|_{\flone } \lesssim  \sum _{k=0}^\infty \frac{2^{k\alpha}C^{k}}{k!} <
\infty.$$ We have proved that the multiplier $e^{i\mu(\xi)} \chi(\xi )$ is in $\flone $ and thus by
Lemma~\ref{suffmulthm} it is bounded on all \modsp s $\M{p,q}$, and clearly also on all Lebesgue
spaces $L^p$.
\end{proof}

\subsection{Large Oscillations at Infinity}

To deal with the oscillatory behavior of $e^{i|\xi|^{\alpha}}$ at infinity, we use a time-frequency
version of the stationary phase method ~\cite{Gra03,Ste93}, the proof of which is reminiscent of
the localization principle for oscillatory integrals of the first kind.  A key feature of the
modulation space case is that we may make a linear alteration of phase without affecting the norm
in $\wam$.

\begin{lemma}\label{lem:ch01}
  Assume that $\alpha (x), \beta (x)$ are arbitrary (measurable)
  functions on $\rd $. Set $\tilde{\sigma}_x (\xi ) = \sigma (\xi )
  e^{i(\alpha (x) + \xi \beta (x))}$. Then
  \begin{equation}
    \label{eq:c3}
    \|\sigma \, T_x g \|_{\flone } = \|\tilde{\sigma}_x \, T_x g
    \|_{\flone } \qquad \forall x\in \rd \, .
  \end{equation}
Consequently, $\|\sigma \|_{\wam } = \sup _{x\in \rd } \|\tilde{\sigma}_x \, T_x g \|_{\flone}$.
\end{lemma}
\begin{proof}
We have
$$ (\tilde{\sigma}_x  \, T_x g )\hat \, (\omega ) = e^{i\alpha (x) }
\intrd \sigma (\xi ) g(\xi-x) e^{-2\pi i \xi (\omega - \frac{\beta (x)}{2\pi }) } \, d\xi =
e^{i\alpha (x) } V_g \sigma \Bigl( x, \omega - \frac{\beta (x)}{2\pi }\Bigr)$$ so $\|\sigma
\|_{\wam } = \sup _x \| \tilde{\sigma}_x \, T_x g \|_{\flone }$ by the translation invariance of
$L^1(d\omega )$.
\end{proof}

Our main result describing the behavior of multipliers with large oscillations is as follows.
\begin{thm}
\label{lemma:c2d} For $d\geq 1$, let $l = \lfloor d/2 \rfloor +1$.  Assume that $\mu $   is
$2l$-times differentiable and $\|\partial ^\beta \mu\|_{L^\infty}\leq C$, for $2\leq |\beta |\leq
2l$, and some  constants $C$. Then $\sigma = e^{i\mu } \in \wam $ and therefore $H_\sigma$ is
bounded on all \modsp s $\M{p,q}$ for $1\leq p, q \leq \infty $.
\end{thm}

\begin{proof}
The argument is most easily understood in the case $d=1$ where the assumption is that $\mu '' \in
L^\infty (\R )$.  We give this proof first and then indicate the necessary modifications for
general $d$.

Let $g$ be a compactly supported test function in $C^\infty $. If we modify the phase $\mu $ by
subtracting the linear Taylor polynomial at $x$
  \begin{equation}
    \label{eq:c5}
    r_x(\xi ) = \mu (\xi ) - \mu (x) - \mu '(x ) (\xi  - x)
  \end{equation}
then by Lemma~\ref{lem:ch01} we have
\begin{align}
    \|\sigma \|_{\wam }
    &= \sup _{x \in \mathbb{R}} \| e^{ir_x} T_xg\|_{\flone} \notag\\
    &= \sup _{x \in \mathbb{R}} \,  \int _{\R } \Big| \int _{\R } e^{i r_x (\xi )}
     T_x g(\xi ) \,e^{-2\pi i \xi \omega } \, d\xi \Big| \, d\omega \notag\\
    &= \sup _{x \in \mathbb{R}}  \int _{|\omega | \leq 1 } \dots \, d\omega +
     \int _{|\omega | \geq 1 } \dots \, d\omega = A + B. \label{eq:decompinoscintegral}
\end{align}

By pulling in the absolute values, the first term is readily estimated by
$$
A \leq \sup _{x \in \mathbb{R}}  \int _{|\omega | \leq 1} \Big( \int _{\R } |g(\xi - x) | d\xi
\Big) \, d\omega \leq 2 \|g \|_{L^1}.
$$

For the estimate of $B$ we write the exponential as $e^{-2\pi i \xi
  \omega } = \frac{-1}{4\pi ^2 \omega ^2} \, \frac{d^2}{d\xi ^2}
e^{-2\pi i \xi \omega } $. Using integration by parts, we obtain
\begin{equation}\label{eq:estforBinosclemma}
    B = \sup _{x \in \mathbb{R}}  \int _{|\omega |\geq 1 } \frac{1}{4\pi ^2 \omega ^2} \Big| \int
  _{\R }\frac{d^2}{d\xi ^2}\Big( e^{i r_x (\xi )} T_x g(\xi )\Big)
\,e^{-2\pi i \xi \omega } \, d\xi \Big| \, d\omega .
\end{equation}
The second derivative in this integral is
\begin{equation*}
\frac{d^{2}}{d\xi^{2}}[T_x g(\xi)\, e^{ir_x( \xi)}]=\bigparen{T_x g^{''}(\xi)  + 2i r_x'(\xi)T_x
g'(\xi) + i r_x^{''}(\xi)T_x g(\xi) - (r_x'(\xi))^{2}T_x g(\xi)}e^{i r_x(\xi)}.
\end{equation*}
However Taylor's theorem supplies bounds $|r_x(\xi)|\lesssim\|\mu''\|_{\infty}|x-\xi|^{2}$ and
$|r_x'(\xi)|\leq \|\mu''\|_{\infty}|x-\xi|$, and it is obvious that
$\|r_x''\|_{\infty}=\|\mu''\|_{\infty}$. Since $T_x g$ and its derivatives are supported in a fixed
neighborhood of $x$ we conclude that
$$
\Big| \int
  _{\R }\frac{d^2}{d\xi ^2}\Big( e^{i r_x (\xi )} T_x g(\xi )\Big)
\,e^{-2\pi i \xi \omega } \, d\xi \Big| \leq C$$ and substituting into \eqref{eq:estforBinosclemma}
we have the bound
$$
B \leq \int _{|\omega | \geq 1} \frac{C}{4\pi ^2 \omega ^2} \, d\omega < \infty \, .
$$
Combining the estimates for $A$ and $B$, we have shown that $e^{i\mu } \in \wam $, and by
Lemma~\ref{suffmulthm} the associated Fourier  multiplier is bounded on $\M{p,q}$ for $1\leq p,q
\leq \infty$.  This concludes the proof for the case $d=1$.

The proof for general $d$ is very similar.  We define
  \begin{equation*}
    r_x(\xi ) = \mu (\xi ) - \mu (x) - \nabla{\mu }(x ) \cdot  (\xi  - x) \, .
  \end{equation*}
and compute as in \eqref{eq:decompinoscintegral}.  Evidently the estimate for the first term
becomes $ A \leq v_d \|g \|_{L^1}$, where $v_d$ is the volume of the unit ball in $\R ^d$. For the
estimate of $B$ we write the exponential as $e^{-2\pi i \xi \cdot  \omega } = \Big(\frac{-1}{4\pi
^2 |\omega| ^2}\Big) ^l \, \Delta ^l ( e^{-2\pi i \xi \cdot \omega })$ and integrate by parts as
before, finding that
\begin{equation*}
    B
    = \sup _{x \in \mathbb{R}^d}  \int _{|\omega |\geq 1 } \frac{1}{(4\pi ^2 |\omega| ^2)^{l}}\,  \Big| \int
  _{\R ^d } \Delta ^l \Big( e^{i r_x (\xi )} T_x g(\xi )\Big)
\,e^{-2\pi i \xi \cdot \omega } \, d\xi \Big| \, d\omega \, .
\end{equation*}
Since $l=\lfloor d/2 \rfloor +1$ we know that $|\omega|^{-2l}$ is integrable outside a neighborhood
of the origin, so the desired bound will follow if we can show  $ \Delta ^l \big( e^{i r_x (\xi )}
T_x g(\xi )\big)$ is uniformly bounded on $\R ^d$.

Using Leibniz's rule and the chain rule, we write $ \Delta ^l \big( e^{i r_x (\xi )} T_x g(\xi
)\big)$ as a linear combination of terms of the form
\begin{equation*}
    e^{ir_x}  \partial ^{\gamma _1} (T_x g ) \, \partial ^{\gamma _2}
     r_x \, \partial ^{\gamma _3} r_x  \dots \partial ^{\gamma _m} r_x
     \qquad   \text{for } \, \sum _{j=1}^m {\gamma _j} = 2l
\end{equation*}
It is immediate for multi-indices $|\gamma _j | \geq 2$ that $|\partial ^{\gamma _j} r_x| =
|\partial ^{\gamma _j} \mu| \leq C$ independent of $x$. Otherwise we may apply Taylor's theorem to
find $|\partial^{\gamma_{j}}r_{x}(\xi)|\leq C|x-\xi|$ when $|\gamma_{j}|=1$ and $|r_x(\xi)|\leq
C|x-\xi|^{2}$. Moreover we are only interested in these functions for $\xi $ in the support of $T_x
g$, on which they are uniformly bounded.  Since the factors $\partial ^{\gamma _1} (T_x g )$ are
also bounded we find that $ \Delta ^l \big( e^{i r_x (\xi )} T_x g(\xi )\big)$ is uniformly
bounded, which gives the desired bound for $B$.

Combining the estimates for $A$ and $B$ we conclude that $e^{i\mu } \in \wam $. By
Lemma~\ref{suffmulthm} the multiplier $H_{e^{i\mu}}$ is then bounded on $\M{p,q}$ for $1\leq p,q
\leq \infty$.
\end{proof}

\subsection{Proof of Theorem \protect{\ref{maintheorem}}}
\begin{proof}
Let $\chi \in \mathcal{C}^{\infty}_{c}(\rd{})$ be the  test function defined in \eqref{bump}. We
split the multiplier into two parts by writing $\sigma(\xi) = e^{i|\xi |^\alpha }\chi(\xi)+
e^{i|\xi |^\alpha }  (1-\chi(\xi)) = \sigma _{\mathrm{sing}}(\xi) + \sigma _{\mathrm{osc}}(\xi).$
Then $\sigma _{\mathrm{sing}} \in \flone $ by Theorem~\ref{lem:ch4} and $H_{\sigma
_{\mathrm{sing}}} $ is bounded on all \modsp s.

To deal with $\sigma _{\mathrm{osc}}$, we set  $\tilde{\mu }(\xi ) = |\xi |^\alpha (1- \chi (2\xi
))$. Then $\tilde{\mu }(\xi ) =  |\xi |^\alpha$ for $|\xi | \geq 1$ and
$$
\sigma _{\mathrm{osc}} (\xi ) = e^{i\tilde{\mu} (\xi) } \big( 1 - \chi (\xi )\big) \, .
$$
By this construction we have removed the singularity of $|\xi |^\alpha $ at the origin.  Since
$\alpha \leq 2$, all derivatives $\partial ^\beta \tilde{\mu }$ are bounded for $|\beta | \geq 2$.
The multiplier $e^{i \tilde{\mu}}$ is therefore bounded on all \modsp s by Theorem~\ref{lemma:c2d}.
Clearly, the multiplier $1-\chi \in \M{\infty ,1}$ is bounded on all \modsp s, so the fact that the
bounded multipliers on $\M{p,q}$ form an algebra implies that the same is true of
$\sigma_{\mathrm{osc}}$. This completes the proof.
\end{proof}

It is not hard to see that we have in fact proven a more general result than
Theorem~\ref{maintheorem}.  Inspecting the conditions needed in Theorems \ref{lem:ch4} and
\ref{lemma:c2d}, we have shown the following.

\begin{cor}\label{metathm}
Let $d\geq 1$, $\alpha \in [0, 2]$, and define $l=\lfloor d/2 \rfloor +1$. Assume that $\mu \in
\mathcal{C}^{2l}(\rd{}\setminus \{0\})$ is homogeneous of order $\alpha$ and all derivatives
$\partial ^\beta \mu $  are bounded outside a neighborhood of $0$ for $2\leq |\beta | \leq 2l$.
Then  $H_{e^{i \mu}}$ is bounded on all modulation spaces $\M{p,q}$ for $1 \leq p, q \leq \infty$.
\end{cor}

For $1< r < \infty$, we let $|\xi|_{r}=(\sum_{j=1}^{d}|\xi_{j}|^{r})^{1/r}$ denote the $r$-norm on
$\rd{}$. With this notation, $|\xi|=|\xi|_{2}$. For  $1\leq r < \infty$, and $\alpha \geq 0$, let
$\mu(\xi)=e^{i |\xi|_{2r}^{\alpha}}.$ Then $\mu(\xi)=|\xi|^{\alpha}_{2r}$ satisfies the conditions
of Corollary \ref{metathm} for all $0\leq \alpha \leq 2 $.
\begin{cor}
The multiplier $e^{i|\xi|_{2r}^{\alpha}}$ is a bounded Fourier multiplier for all modulation
spaces.
  \end{cor}

\subsection{Improved Estimates  for  the Cases
  $\alpha=2$ and $\alpha\in [0, 1]$}

Boundedness of the Fourier multiplier operator $H_{\sigma_{2}}$ on $L^{p}(\rd{})$ was settled by
H\"ormander \cite{Ho63}, who showed the more general result that when $\phi$ is any quadratic
polynomial the multiplier $\sigma = e^{i\phi}$ is bounded only on $L^{2}(\rd{})$. We now give a
different proof of part of Theorem \ref{maintheorem}, showing boundedness of $H_{\sigma _2}$ by a
time-frequency approach that is based on the so-called metaplectic invariance of the modulation
spaces.  This method gives a better bound for the operator norm.

\begin{thm}\label{alpha2}
Let $d\geq 1$, and let $\sigma_2(\xi) = e^{i\pi t|\xi|^{2}}$. Then:

\noindent (a) $\sigma_2\in (\wam \cup \M{1, \infty} )\setminus \M{\infty, 1}$,

\noindent (b) $H_{\sigma_2}$ is bounded on all modulation spaces $\M{p,q}(\rd{})$, $1\leq p, q \leq
\infty$, and the operator norm satisfies the uniform estimate
$\nm{H_{\sigma_{2}}}{{\textrm{op}}}\leq c(d, p, q) (1 + t^2)^{d/4}$.

\noindent (c) Let  $\sigma (\xi ) =  e^{-\pi \xi \cdot A\xi +2\pi b\cdot
  \xi}$ be  a generalized Gaussian so that $A=B+iC $ for a
positive-definite real-valued $d\times d$-matrix $B$, a symmetric real-valued matrix $C$ and $b \in
\mathbb{C} ^d$.  Then $H_\sigma $ is  bounded on all modulation spaces $\M{p, q}(\rd{})$, $1\leq p,
q \leq \infty$.
\end{thm}

\begin{proof}
We use the Gaussian  $g(\xi) = e^{-\pi |\xi|^{2}}$ as a window for the short-time Fourier
transform. Then the STFT $V_g \sigma $ can be calculated explicitly by using  Gaussian integrals.
\begin{align*}
V_{g}\sigma_{2}(x, \omega) & = \int_{\rd{}}e^{i\pi t |\xi|^2}\, e^{-2\pi i \xi \cdot \omega}\, e^{-\pi|\xi - x|^2}\, d\xi\\
& = e^{-\pi |x|^2}\, \int_{\rd{}}e^{-\pi (1-it) |\xi|^{2}} e^{  2\pi
  \xi \cdot x} \, e^{-2\pi i \xi \cdot \omega}\, d\omega \, .
\end{align*}
The integral is the Fourier transform of a generalized Gaussian. By using a table or
\cite[Lemma~4.4.2]{Gr01}, we obtain
\begin{eqnarray*}
  V_{g}\sigma_{2}(x, \omega) & = & e^{-\pi |x|^2}\, (1-it)^{-d/2} \,
  e^{\pi (1-it) |x|^2} \, T_{tx} M_{-x} \Big( e^{-\pi |\omega |^2/
    (1-it)} \Big) \,
\end{eqnarray*}
(where the square root $(1-it)^{1/2}$ is taken with positive imaginary part). After taking absolute
values and performing some cancellations we arrive at the expression
\begin{equation}
  \label{eq:ch13}
  |V_{g}\sigma_{2}(x, \omega)| = (1+t^2)^{-d/4} e^{-\pi |\omega
    -tx|^2/(1+t^2)} \, .
\end{equation}
Since $\int _{\rd } e^{-a|x|^2}\,dx= a^{-d/2}$, the modulation space norms of $\sigma _2$ are now
easy to compute.  It is trivial that
$$\int_{\rd{}}\nm{V_{g}\sigma_{2}(\cdot, \omega)}{L^{\infty}}\,
d\omega =  \infty,$$ and therefore $\sigma_2 \not\in \M{\infty, 1}(\rd{})$. On the other hand,
$$ \|\sigma _2\|_{\M{1,\infty }} = \sup _\omega
\int_{\rd{}}|V_{g}\sigma_{2}(x, \omega)|\, dx = (1+t^2)^{-d/4} \int_{\rd{}}
 e^{-\frac{\pi t^2}{t^2 +1}|x|^{2}}\, dx = (t^2 +1)^{d/4} t^{-d},$$
and
$$
 \|\sigma _2\|_{\wam } = \sup _x
\int_{\rd{}}|V_{g}\sigma_{2}(x, \omega)|\, d\omega = (1+t^2)^{-d/4} \int_{\rd{}}
 e^{-\frac{\pi }{t^2 +1}|\omega|^{2}}\, d\omega = (t^2 +1)^{d/4} \, .$$
Consequently $\sigma _2 \in \M{1,\infty}$ and $\sigma _2 \in \wam $, so Lemma~\ref{suffmulthm}
implies the boundedness of $\sigma _2$ with an explicit form for the dependence on the parameter
$t$:
$$
\|H_{\sigma _2} f\|_{\M{p,q}} \lesssim \|\sigma _2\|_{\wam } \|f\|_{\M{p,q}} \lesssim (1+t^2)^{d/4}
\|f\|_{\M{p,q}} \, .
$$

The proof of (c) is similar, using the fact that after a change of coordinates, any quadratic
function on $\R^d$ can be written in the form $\phi(\xi) = \ip{\xi}{C\xi}$, where $C$ is a $d
\times d$ hermitian matrix.
\end{proof}

For the range $\alpha \in [0,1]$, we now prove a stronger property of the multipliers $e^{i|\xi
|^\alpha }$, which is also sufficient to show that $H_{\sigma _\alpha } $ is bounded on all
modulation spaces.

\begin{cor}\label{alpha01}
If  $\alpha \in [0, 1]$, then $\sigma_{\alpha}(\xi) =e^{i|\xi|^{\alpha}}$ belongs to $ \M{\infty,
  1}(\rd{})$.
\end{cor}

\begin{proof} Let $\chi$ be the smooth bump function defined by ~\eqref{bump}, and write
$$\sigma_{\alpha}(\xi)
= \chi(\xi) \sigma_{\alpha}(\xi) + (1-\chi(\xi))\sigma_{\alpha}(\xi) = \sigma_{\mathrm{sing}}(\xi)
+ \sigma_{\mathrm{osc}}(\xi).$$ Theorem \ref{lem:ch4} implies that $\sigma_{\mathrm{sing}} \in
\flone \subset \M{\infty, 1}$. It is readily seen that  the following estimate holds  for all
$|\xi|\geq 1$ and $|\beta| \geq 1$:
$$|\partial^{\beta}\sigma_{\mathrm{osc}}(\xi)|\leq C_{\beta}(1 +
|\xi|)^{\alpha - |\beta|} \, .$$
 Since $\alpha \leq 1$, all partial derivatives of $\sigma
 _{\mathrm{osc}}$ are bounded, and this fact implies that
$\sigma_{\mathrm{osc}} \in \mathcal{C}^{d+1}(\rd{}) \subset \M{\infty, 1}(\rd{})$. For this
embedding, see, e.g., \cite[Thm.~14.5.3]{Gr01} or \cite{Ok04}.  Thus $\sigma _\alpha \in \M{\infty
,1}$ and the conclusion  follows.
\end{proof}

\subsection{Further results}

Next  we  consider  the related  family of multipliers $\sigma_{\alpha, \delta}$ defined by
$$\sigma_{\alpha, \delta}(\xi)=\Im (e^{i|\xi|^\alpha}|\xi|^{-\delta})=\frac{\sin|\xi|^\alpha}{|\xi|^{\delta}},
\alpha, \delta>0.$$ The following statement  should be compared to results proved in \cite{Hir,
Per, Wa65}.

\begin{thm}\label{extra}
Let $d\geq 1$, and let $\alpha, \beta >0$.

\noindent (a.) If $0<\delta \leq \alpha\leq 1 $. Then $\sigma_{\alpha, \delta}\in\M{\infty, 1}$.
Consequently, $H_{\sigma_{\alpha, \delta}}$ is bounded on $M^{p, q}$ for all $1\leq p, q\leq
\infty$

 \noindent (b) If $\alpha >1$ and $\delta \leq \alpha$, then
$H_{\sigma_{\alpha, \delta}}$ is bounded on $M^{p, q}$ for $1\leq p, q\leq \infty$ and
$|\frac{1}{p}-\frac{1}{2}|< \frac{\delta}{\alpha d}.$
\end{thm}

\begin{proof} (a) Using the smooth bump $\chi $  defined by ~\eqref{bump},
  we  write
$$\sigma_{\alpha, \delta}(\xi)= \chi(\xi) \sigma_{\alpha, \delta}(\xi) + (1-\chi(\xi))\sigma_{\alpha, \delta}(\xi)
= \sigma_{\mathrm{sing}}(\xi) + \sigma_{\mathrm{osc}}(\xi).$$ We first show that
$\sigma_{\mathrm{sing}}\in \flone$. Using the same notation as in Theorem \ref{lem:ch4}, we
decompose the symbol as
$$\sigma_{\mathrm{sing}}=\displaystyle\sum_{k=0}^\infty\frac{(-1)^k}{(2k+1)!}
|\xi |^{(2k+1)\alpha -\delta } \chi (\xi ) \, .$$ Set $\psi_{k}(\xi)=|\xi|^{(2k+1)\alpha-\delta}
\psi (\xi)$ and
$$\phi_{k}=|\xi|^{(2k+1)\alpha-\delta}\chi
(\xi)=\displaystyle\sum_{j=1}^\infty 2^{-j(2k\alpha+\alpha-\delta)}\psi_{k}(2^j\xi)\, .$$ Since
$(2k+1)\alpha -\delta >0$, the proof of Theorem~\ref{lem:ch4} applies and we conclude that $\sigma
_{\mathrm{sing}} \in \flone$ and hence $H_{\sigma
  _{\mathrm{sing}}}$ is bounded on all $\M{p,q}$.

Now consider $\sigma _{\mathrm{osc}}= (\sin |\xi |^\alpha) | \, \xi |^{-\delta } (1-\chi (\xi ))$.
The multiplier $\frac{1}{2i} (e^{i|\xi
  |^\alpha } - e^{-i|\xi
  |^\alpha }  )$ is bounded on all modulation spaces by
Theorem~\ref{maintheorem}. On the other hand, after removing the singularity at $\xi = 0$, the
multiplier $\kappa (\xi ) = |\xi |^{-\delta } (1-\chi (\xi ))$ satisfies the conditions of the
H\"ormander-Mihlin multiplier theorem. In particular, all partial derivatives $\partial ^\beta
\kappa $ are bounded. As before we conclude that $  |\xi |^{-\delta } (1-\chi (\xi )) \in
C^{d+1}(\rd ) \subseteq \M{\infty
  ,1}$, and thus $H_\kappa $ is bounded on all modulation spaces
$\M{p,q}$. Consequently $H_{\sigma _{\mathrm{osc}}} = H_{\sin |\xi
    |^\alpha } H_{\kappa }$ is also bounded on $\M{p,q}$ for $1\leq
  p,q \leq \infty $, and the theorem is proved.

  (b) We still write $\sigma_{\alpha, \delta}(\xi)= \sigma_{\mathrm{sing}}(\xi) +
  \sigma_{\mathrm{osc}}(\xi)$, and the same argument as above
  shows that  $\sigma
_{\mathrm{sing}} \in \flone$ and hence $H_{\sigma
  _{\mathrm{sing}}}$ is bounded on all $\M{p,q}$ (in fact $ H_{\sigma
  _{\mathrm{sing}}}$ is bounded on all $L^{p}(\R^{d})$, $1\leq p
  \leq \infty$). On the other hand it was proved in \cite{Hir, Per,
  Wa65} that $H_{\sigma _{\mathrm{osc}}}$ is bounded on
  $L^{p}(\R^{d})$ whenever $|\frac{1}{p}-\frac{1}{2}|<
\frac{\delta}{\alpha d}.$ Consequently, using \cite[Theorem 17]{FeiNa} we conclude that $H_{\sigma
_{\mathrm{osc}}}$ is bounded on $\M{p,q}$ whenever $|\frac{1}{p}-\frac{1}{2}|< \frac{\delta}{\alpha
d}$ and for all $1\leq q \leq \infty$. This concludes the proof.
\end{proof}
\begin{rem}
In contrast to Theorem~\ref{extra}, the multipliers $\tilde\sigma_{\alpha,
\delta}(\xi)=e^{i\xi|^\alpha}|\xi|^{-\delta}$ for $\alpha ,\delta >0$ and  $|\xi |^{-\delta}\sin
|\xi |^\alpha $ for  $\delta >\alpha >0$ are not bounded on $L^p$ or on $\M{p,q}$, because they are
unbounded functions. Using arguments of this section, we can show that
 the Fourier multiplier with
 symbol $$\tilde\sigma_{\alpha, \delta} - \displaystyle\sum_{k=0}^{\lfloor \delta/\alpha \rfloor}
\frac{i^k\phi_{k, \alpha, \delta}}{k!}$$ is bounded on certain modulation spaces.
\end{rem}

\section{Applications to Some Cauchy Problems}

\subsection{The Schr\"odinger equation}

Consider the linear free Schr\"odinger equation
\begin{equation}
\label{schro} \left\{\begin{array}{r@{\quad =\quad}l}
i\frac{\partial{u}}{\partial{t}}(x,t) & \Delta_{x}u(x, t)\\
u(x, 0) & f(x), x\in\rd, t\geq 0, \end{array}\right.
\end{equation}
where $\Delta_x$ is the  Laplacian.  The formal solution to this equation is given by
\begin{equation}\label{solschro}
u(x, t) = \int_{\rd{}}e^{i t |\xi|^{2}}\, \hat{f}(\xi)\, e^{2 \pi i \xi \cdot x}\,
d\xi=H_{\sigma_{2}^{t}}f(x),
\end{equation} where $\sigma_{2}^{t}(\xi) = \sigma_{2}(\sqrt{t}\xi) e^{i t |\xi|^{2}}.$ is a
bounded multiplier on modulation spaces by Theorem~\ref{maintheorem} and Theorem~\ref{alpha2}.

\begin{cor}\label{estsolsch}
Let $d\geq 1$, and let $u(x, t)$ be given by ~\eqref{solschro}. Then, for any $t\geq 0$,
$$\nm{u(\cdot, t)}{\M{p,q}} \leq C(t^{2} + 4 \pi^{2})^{d/4}\,
\nm{f}{\M{p,q}}$$ for all $1\leq p, q \leq \infty$ and a  constant $C$ depending only on $d, p$ and
$q$.
\end{cor}

\begin{rem}
  This statement was also obtained with a different
  method in ~\cite{wang06}.
\end{rem}
\begin{rem}
In particular, modulation space properties are preserved by the time evolution of  the
Schr\"odinger equation. This is in strong contrast to the standard $L^p$-theory where  the
$L^p$-property of the initial data is not preserved by the time evolution. see, for example,
\cite{Str77}, where it was shown that $L^2(\rd)\ni f(x) \mapsto u (x, t)\in L^p(\R^{d+1})$ for
$p=2(d+2)/d$.
\end{rem}

\subsection{The Wave Equation}

Consider now the following Cauchy problem for the wave equation
\begin{equation}\label{wave}
\left\{\begin{array}{r@{\quad = \quad}l}
\frac{\partial^{2}{u}}{\partial{t^{2}}}(x,t) & \Delta_{x}u(x, t)\\
u(x, 0) & f(x) \\
\frac{\partial{u}}{\partial{t}}(x, 0) & g(x).
\end{array}\right.
\end{equation}
Its formal  solution is  given by
\begin{equation}\label{solwave}
u(x, t) = \int_{\rd{}}\cos{(t |\xi|)}\, \hat{f}(\xi)\, e^{2 \pi i \xi \cdot x}\, d\xi +
 \int_{\rd{}}\frac{\sin{t |\xi|}}{|\xi|}\, \hat{g}(\xi) \, e^{2 \pi i \xi \cdot x}\, d\xi.
 \end{equation}

The time evolution requires an understanding of the continuity properties of the Fourier
multipliers $\sigma^t(\xi) = \cos{t |\xi|}$, or equivalently $\sigma_{1}^t(\xi) = e^{i t |\xi|}$,
and $m^t (\xi) = \frac{\sin{t |\xi|}}{|\xi|}.$ The first of these multipliers is known to be
bounded on all $L^{p}(\R)$, but only on $L^{2}(\rd{})$ for all $d\geq 1$, \cite{Ho63, Lit63}.
Theorems \ref{maintheorem} and~\ref{extra} yield  the following result.

\begin{cor}\label{estsolwave}
Let $d\geq 1$, and let $u(x, t)$ be the solution of the wave equation as  given by
~\eqref{solwave}. Then, for any $t\geq 0$,
$$\nm{u(\cdot, t)}{\M{p,q}} \leq C(t)\, (\nm{f}{\M{p,q}} +
\nm{g}{\M{p, q}})$$ for all $1\leq p, q \leq \infty$, where $C(t)>0$ depends on $ d, p,$ and $ q.$
Again, the solution to the Cauchy problem for the wave equation preserves the  initial data in a
modulation space.
\end{cor}

\begin{rem}
Again, one should compare the space preserving estimate in the previous theorem to, for example,
the following boundedness result $\dot{L_1^p} (\rd{})\times L^p (\rd{}) \ni (f, g) \mapsto u(\cdot,
t) \in L^q (\rd{}),$ for certain $p\leq 2\leq q,$ proved in \cite{Str70a}. Here, $\dot{L_1^p}$
denotes the appropriate homogeneous Sobolev space.
\end{rem}

As previously mentioned, these results show that solutions to the Cauchy problems for the
Schr\"odinger and the wave equation with initial data in a modulation space stay in the same space
for all future time. This is a time-frequency version of the classical principle of conservation of
energy \cite{Lit63, Str70a, Str70b, Str77} for these Cauchy problems. The reader will recall that,
in the context of Lebesgue spaces, this principle holds only on $L^2$.

\section{Acknowledgments}

The authors would like to thank Carlos Kenig  and Robert Strichartz for bringing  some of the
questions discussed in this work to their attention. They also thank Hans Feichtinger and Camil
Muscalu for very helpful discussions. This work was partially developed while the authors were
visiting the Erwin Schr\"odinger Institute (ESI) in Vienna. Its  support and  hospitality are
gratefully acknowledged.


\begin{thebibliography}{20}

\bibitem{wang06}
W.~Baoxiang, Z.~Lifeng, and G.~Boling,
\newblock Isometric decomposition operators, function spaces {$E\sp \lambda\sb
  {p,q}$} and applications to nonlinear evolution equations,
\newblock  J. Funct. Anal., 233(1):1--39, 2006.


\bibitem{Be82}
R.~M.~Beals,
\newblock $L^p$ boundedness of Fourier integral operators,
\newblock Mem.\ Amer.\ Math.\ Soc.\ {\bf 264} (1982).

\bibitem{BGGO}
{\'A}.~B\'enyi, L.~Grafakos, K.~Gr\"ochenig, and K.~Okoudjou,
\newblock A class of {F}ourier multipliers for modulation spaces,
\newblock Appl.\ Comput.\ Harmon.\ Anal.\ {\bf 19} (2005), no.\ 1, 131--139.

\bibitem{CG03}
E.~Cordero and K.~Gr{\"o}chenig,
\newblock Time-frequency analysis of localization operators,
\newblock {\em J. Funct. Anal.}, 205(1):107--131, 2003.

\bibitem{Fei83}
H.~G.~Feichtinger,
\newblock  Modulation spaces on locally {A}belian groups, Technical Report,
University of Vienna, 1983,
\newblock Updated version appeared in Proceedings of ``International Conference on
Wavelets and Applications'' 2002, pp.~99-140, Chennai, India, 2003.

\bibitem{fei90}
H.~ G.~Feichtinger,
\newblock Generalized amalgams, with applications to Fourier transform,
\newblock Can.\ J.\ Math.\ {\bf 42} (1990), 395--409.

\bibitem{fei06}
H.~Feichtinger,
\newblock Modulation spaces: looking back and ahead,
\newblock Sampl. Theory Signal Image Process., 5(2):109--140, 2006.

\bibitem{FG89a}
H.~G.~Feichtinger and K.~Gr\"ochenig,
\newblock  Banach spaces related to integrable group representations and their atomic decompositions I,
\newblock J.~ Funct.\ Anal.\ {\bf 86} (1989), 307--334.

\bibitem{FG89b}
H.~ G.~Feichtinger and K.~Gr\"ochenig,
\newblock Banach spaces related to integrable group representations and their atomic decompositions II,
\newblock Monatsh.\ Math.\ {\bf 108} (1989), 129--148.

\bibitem{FeiNa}
H.~G.~Feichtinger and G.~Narimani,
\newblock Fourier multipliers of classical modulation spaces,
\newblock Appl.\ Comput.\ Harmon.\ Anal., to appear.

\bibitem{Fol89}
G.~B.~Folland,
\newblock ``Harmonic Analysis in Phase Space'',
\newblock Ann.\ of Math.\ Studies, Princeton University Press, Princeton NJ, 1989.

\bibitem{Gra03}
L.~Grafakos,
\newblock  ``Classical and Modern Fourier Analysis'',
\newblock Prentice Hall, Upper Saddle River NJ, 2003.

\bibitem{Gr01}
K.~Gr\"ochenig,
\newblock ``Foundations of Time-Frequency Analysis'',
\newblock Birkh\"auser, Boston MA, 2001.

\bibitem{grhe}
K.~Gr\"ochenig and C.~Heil,
\newblock Modulation spaces and pseudodifferential operators
\newblock, Int.\ Eq.\ Oper.\ Theory \textbf{34} (1999), 439--457.

\bibitem{gro06}
K.~Gr\"ochenig,
\newblock Time-frequency analysis of {S}j\"ostrand's class,
\newblock {\em Revista Mat. Iberoam.}, 22(2):703--724, 2006,
\newblock arXiv:math.FA/0409280v1.

\bibitem{Hir}
I.~I.~ Hirschman,
\newblock On multiplier transformations,
\newblock Duke Math.\ J.\ {\bf 26} (1959), 221--242.

\bibitem{Ho63}
L.~H\"ormander,
\newblock Estimates for translation invariant operators in $L^{p}$ spaces,
\newblock Acta Math.\ {\bf 104} (1960), 93--140.

\bibitem{LeOl}
V.~Lebedev and A.~Olevski\v{i},
\newblock $\mathcal{C}^{1}$ changes of variable: Beurling-Helson type theorem and H\"ormander conjecture on Fourier multipliers,
\newblock Geom.\ Funct.\ Anal.\ {\bf 4} (1994), no.\ 2, 213--235

\bibitem{Lit63}
W.~Littman,
\newblock The wave operator and $L_{p}$ norms,
\newblock J.\ Math.\ Mech.\ {\bf 12} (1963), 55--68.

\bibitem{mi56}
S.~G.~Mihlin,
\newblock On the multipliers of Fourier integrals,
\newblock Dokl.\ Akad.\ Nauk.\ SSSR (N.\ S.\ ),{\bf 109} (1956), 701--703 (Russian).

\bibitem{Ok04}
K.~A.~Okoudjou,
\newblock Embeddings of some classical Banach spaces into the modulation spaces,
\newblock  Proc.\ Amer.\ Math.\ Soc., 132 (2004), no.\ 6, 1639--1647.

\bibitem{Per}
J.~C.~Peral,
\newblock $L^p$ estimates for the wave equation,
\newblock J.\ Funct.\ Anal.\ {\bf 36} (1980), 114--145.

\bibitem{Sj94}
J.~Sj\"ostrand,
\newblock An algebra of pseudodifferential operators,
\newblock Math.\ Res.\ Lett.\ {\bf 1} (1994), 185--192.

\bibitem{Sog}
C.~D.~Sogge,
\newblock $L^p$ estimates for the wave equation and applications,
\newblock  Journ\'ees "\'Equations aux D\'eriv\'ees Partielles" (Saint-Jean-de-Monts, 1993),  Exp. No. XV, 1--12, \'Ecole Polytech., Palaiseau, 1993.

\bibitem{Ste70}
E.~M.~Stein,
\newblock ``Singular integrals and differentiability properties of functions'',
\newblock  Princeton Mathematical Series, Princeton University Press, Princeton, N.J. 1970.

\bibitem{Ste93}
E.~M.~Stein,
\newblock ``Harmonic Analysis: Real-Variable Methods, Orthogonality, and Oscillatory Integrals'',
\newblock Princeton University Press, Princeton, NJ, 1993.

\bibitem{Str70a}
R.~Strichartz,
\newblock Convolutions with kernels having singularities on a sphere,
\newblock Trans.\ Amer.\ Math.\ Soc.\ {\bf 148} (1970), 461--471.

\bibitem{Str70b}
R.~S.~Strichartz,
\newblock A priori estimates for the wave equation and some applications,
\newblock J.\ Funct.\  Anal.\  {\bf 5} (1970), 218--235.

\bibitem{Str77}
R.~S.~Strichartz,
\newblock Restrictions of Fourier transforms to quadratic surfaces and decay of solutions of wave equations,
\newblock Duke Math.\ J.\  {\bf 44} (1977), no.\ 3, 705--714.

\bibitem{toft}
J.~Toft,
\newblock Continuity properties for modulation spaces, with applications to pseudodifferential operators. I,
\newblock J.\ Funct.\ Anal.\ 207 (2004), 399--429.

\bibitem{Wa65}
S.~Wainger,
\newblock Special trigonometric series in $k$-dimensions,
\newblock Mem.\ Amer.\ Math.\ Soc.\ {\bf 59} (1965).

\end{thebibliography}
\end{document}